\documentclass[12pt]{report}
\usepackage{amssymb,amsmath,makeidx,verbatim}
\newcommand{\R}{\mathbb{R}}

\newcommand{\C}{\mathbb{C}}
\newcommand{\Z}{\mathbb{Z}}

\newcommand{\be}{\begin{enumerate}}
\newcommand{\ee}{\end{enumerate}}
\newcommand{\bq}{\begin{eqnarray*}}
\newcommand{\eq}{\end{eqnarray*}}

\begin{document}
\newcommand{\disp}{\displaystyle}
\thispagestyle{empty}
\begin{center}
\textsc{Functional analysis of canonical wave-packets on real reductive groups\\}
\ \\
\textsc{Olufemi O. Oyadare}\\
\ \\
Department of Mathematics,\\
Obafemi Awolowo University,\\
Ile-Ife, $220005,$ NIGERIA.\\
\text{E-mail: \textit{femi\_oya@yahoo.com}}\\
\end{center}
\begin{quote}
{\bf Abstract.} {\it This paper revisits the study of canonical wave-packets (on a real reductive group $G$ in the Harish-Chandra class) constructed by Oyadare $[9c.]$ for the proof of the fundamental theorem of harmonic analysis on $G$ given, in terms of the convolutions $\ast$ on $G$ and the Harish-Chandra transform $\mathcal{H}$ on $\mathcal{C}^{p}(G),$ as $$\psi_{a}=\mathcal{H}^{-1}((\mathcal{H}\xi_{1})^{-1})\ast\mathcal{H}^{-1}((\mathcal{H}\xi_{1})\cdot a\cdot(\mathcal{H}\xi_{1}))\ast\mathcal{H}^{-1}((\mathcal{H}\xi_{1})^{-1})$$ with $a\in\mathcal{C}^{p}(\widehat{G}).$ Following an extraction of its spherical and symmetric convolution-components, a simpler expression is derived for $\psi_{a}$ from which some of its properties are established. Among other results, we employ the wave-packets to give decompositions of the Schwartz-type algebras $\mathcal{C}^{p}(G).$}
\end{quote}

\ \\
\ \\
\ \\
\ \\
\ \\
\ \\
\ \\
\ \\
\ \\
\ \\
\ \\
$\overline{2010\; \textmd{Mathematics}}$ Subject Classification: $43A85, \;\; 22E30, \;\; 22E46$

{\bf \S1. Introduction.}

Let $G$ be a reductive group in the \textit{Harish-Chandra class} where $\mathcal{C}^{p}(G)$ is the \textit{Harish-Chandra-type} Schwartz algebra on $G,$ $0<p\leq2,$ with $\mathcal{C}^{2}(G)=:\mathcal{C}(G).$ It is known that $C^{\infty}_{c}(G)$ is dense in $\mathcal{C}^{p}(G),$ with continuous inclusion. The image of $\mathcal{C}^{p}(G)$ under the (Harish-Chandra) \textit{Fourier transform} on $G$ has been a pre-occupation of harmonic analysts since Harish-Chandra defined $\mathcal{C}(G)$ leading to the emergence of Arthur's thesis $[1a],$ where the Fourier image of $\mathcal{C}(G)$ was characterized for connected non-compact semisimple Lie groups of real rank one. Thereafter Eguchi $[3a.]$ removed the restriction of the real rank and considered non-compact real semisimple $G$ with only one conjugacy class of \textit{Cartan subgroups} as well as the Fourier image of $\mathcal{C}^{p}(G/K)$ in $[3b.],$ while Barker $[2.]$ considered $\mathcal{C}^{p}(SL(2,\R))$ as well as the zero-Schwartz space $\mathcal{C}^{0}(SL(2,\R)).$

The complete $p=2$ story for any real reductive $G$ is contained in Arthur $[1b,c].$ The most successful general result along the general case of $p$ is the well-known \textit{Trombi-Varadarajan Theorem} $[11.]$ which characterized the image of $\mathcal{C}^{p}(G//K),$ $0<p\leq2,$ for a maximal compact subgroup $K$ of a connected semisimple Lie group $G$ as a (Schwartz) multiplication algebra $\bar{\mathcal{Z}}({\mathfrak{F}}^{\epsilon})$ (of $\mathfrak{w}-$invariant members of $\mathcal{Z}({\mathfrak{F}}^{\epsilon}),$ with $\epsilon=(2/p)-1$); thus subsuming the works of Ehrepreis and Mautner $[5.]$ and Helgason $[7.].$ However the characterization of the image of $\mathcal{C}^{p}(G)$ for reductive groups $G$ in the Harish-Chandra class had been long sought, due to failure of the method of generalizing from the real rank one case (successfully employed in $[1b.,c.],$ $[3a.]$ and $[10c.]$) or from the spherical case (considered in $[11.]$). This characterization, also called \textit{the fundamental theorem of harmonic analysis on $G,$} has been recently announced in Oyadare $ [9c.],$ thus paving the way for the proof of the full Bochner theorem on $G,$ $[9d.].$ This paper contains a more detailed considerations of the properties of the wave-packets of the said fundamental theorem.

\ \\
{\bf \S2. Structure of groups in the Harish-Chandra class}

Let $G$ be a group in the Harish-Chandra class. That is $G$ is a locally compact group with the properties that $G$ is reductive, with Lie algebra $\mathfrak{g},$ $[G:G^{0}]<\infty,$ where $G^{0}$ is the connected component of $G$ containing the identity, in which the analytic subgroup, $G_{1},$ of $G$ defined by $\mathfrak{g}_{1}=[\mathfrak{g},\mathfrak{g}]$ is closed in $G$ and of finite center and in which, if $G_{\C}$ is the adjoint group of $\mathfrak{g}_{\C},$ then $Ad(G)\subset G_{\C}.$ Such a group $G$ is endowed with a \textit{Cartan involution,} $\theta,$ whose fixed points form a \textit{maximal} compact subgroup, $K,$ of $G$ $[6.].$ $K$ meets all connected components of $G,$ in particular $K\cap G^{0}\neq\phi.$ Let $\mathfrak{t}$ denote the Lie algebra of $K.$

We denote the \textit{universal enveloping algebra} of $\mathfrak{g}_{\C}$ by $\mathcal{U}(\mathfrak{g}_{\C}),$ whose members may be viewed either as left or right invariant differential operators on $G.$ We shall write $f(x;a)$ for the left action $(af)(x)$ and $f(a;x)$ for the right action $(fa)(x)$ of $\mathcal{U}(\mathfrak{g}_{\C})$ on functions $f$ on $G.$ Let $\mathcal{C}(G)$ represents the space of $C^{\infty}-$functions $f$ on $G$ for which $$\sup_{x \in G}\mid f(b;x;a) \mid \Xi^{-1}(x)(1+\sigma(x))^{r}<\infty,$$ for $a,b \in \mathcal{U}(\mathfrak{g}_{\C})$ and $r>0.$ Here $\Xi$ and $\sigma$ are well-known \textit{elementary} spherical functions defined below on $G.$ $\mathcal{C}(G)$ is known to be a Schwartz algebra under convolution while $\mathcal{C}(G//K),$ consisting of the spherical members of $\mathcal{C}(G),$ is a closed commutative subalgebra. $C^{\infty}_{c}(G)$ is densely contained in $\mathcal{C}(G),$ with continuous inclusion.

Let $\hat{G}$ represent the set of equivalence classes of \textit{irreducible unitary representations} of $G.$ If $G_{1}$ is non-compact then the support of the Plancherel measure does not exhaust $\hat{G}.$ We write $\hat{G_{t}}$ for this support, which generally contains a discrete part, $\hat{G_{d}}$ ($\neq\emptyset,$ if $rank(G) = rank(K)$), and a continuous part, $\hat{G_{t}}\setminus \hat{G_{d}}$ ($\neq\emptyset,$ always).

If $\mathfrak{p}=\{X\in\mathfrak{g}:\theta X=-X\}$ then $\mathfrak{g}=\mathfrak{t}\oplus\mathfrak{p}.$ Choose a maximal abelian subspace  $\mathfrak{a}$ of $\mathfrak{p}$ with algebraic
dual $\mathfrak{a}^*$ and set $A =\exp \mathfrak{a}.$  For every $\lambda \in \mathfrak{a}^*$ put
$$\disp\mathfrak{g}_{\lambda} = \{X \in \mathfrak{g}: [H, X] =
\lambda(H)X, \forall  H \in \mathfrak{a}\},$$ and call $\lambda$ a \textit{restricted
root} of $(\mathfrak{g},\mathfrak{a})$ whenever $\mathfrak{g}_{\lambda}\neq\{0\}$.
Denote by $\mathfrak{a}'$ the open subset of $\mathfrak{a}$
where all restricted roots are $\neq 0$,  and call its connected
components the \textit{Weyl chambers}.  Let $\mathfrak{a}^+$ be one of the Weyl
chambers, define the restricted root $\lambda$ positive whenever it
is positive on $\mathfrak{a}^+$ and denote by $\triangle^+$ the set of all
restricted positive roots.  We then have the \textit{Iwasawa
decomposition} $G = KAN$, where $N$ is the analytic subgroup of $G$
corresponding to $\disp \mathfrak{n} = \sum_{\lambda \in \triangle^+} \mathfrak{g}_{\lambda},$
and the \textit{polar decomposition} $G = K\cdot
cl(A^+)\cdot K,$ with $A^+ = \exp \mathfrak{a}^+,$ and $cl(A^{+})$ denoting the closure of $A^{+}.$

If we set $\disp M = \{k
\in K: Ad(k)H = H, H\in \mathfrak{a}\}$ and $\disp M' = \{k
\in K : Ad(k)\mathfrak{a} \subset \mathfrak{a}\}$ and call them the
\textit{centralizer} and \textit{normalizer} of $\mathfrak{a}$ in $K,$ respectively, then;
(i) $M$ and $M'$ are compact and have the same Lie algebra and
(ii) the factor  $\mathfrak{w} = M'/M$ is a finite group called the \textit{Weyl
group}.  $\mathfrak{w}$ acts on $\mathfrak{a}^*_{\C}$ as a group of linear
transformations by the requirement $$(s\lambda)(H) =
\lambda(s^{-1}H),$$ $H \in \mathfrak{a}$, $s \in \mathfrak{w}$, $\lambda \in
\mathfrak{a}^*_\mathbb{\C}$, the complexification of $\mathfrak{a}^*$.  We then have the
\textit{Bruhat decomposition} $$\disp G = \bigsqcup_{s\in \mathfrak{w}} (B m_sB)$$ where
$B = MAN$ is a closed subgroup of $G$ and $m_s \in M'$ is the
representative of $s$ (i.e., $s = m_sM$).

Some of the most important functions on $G$ are the \textit{spherical
functions} which we now discuss as follows.  A non-zero continuous
function $\varphi$ on $G$ shall be called \textit{(elementary or zonal) spherical
function} whenever

$(i.)\;\varphi(e)=1,$ $$(ii.)\;\varphi \in C(G//K):=\{g\in
C(G): g(k_1 x k_2) = g(x), k_1,k_2 \in K, x \in G\}$$ and $(iii.)\;f*\varphi
= (f*\varphi)(e)\cdot \varphi$ for every $f \in C_c(G//K)$.  This
leads to the existence of a homomorphism $\lambda :
C_c(G//K)\rightarrow \C$ given as $\lambda(f) = (f*\varphi)(e)$.
This definition of an elementary spherical function is equivalent to the functional relation $$\disp
\int_K\varphi(xky)dk = \varphi(x)\varphi(y),$$ $x,y\in G$.  It has
been shown by Harish-Chandra [$6.$] that elementary spherical functions on $G$
can be parametrized by members of $\mathfrak{a}^*_{\C}$.  Indeed every elementary
spherical function on $G$ is of the form $$\disp
\varphi_{\lambda}(x) = \int_Ke^{(i\lambda-p)H(xk)}dk,\; \lambda
\in \mathfrak{a}^*_{\C},$$  $\disp \rho =
\frac{1}{2}\sum_{\lambda\in\triangle^+} m_{\lambda}\cdot\lambda,$ where
$m_{\lambda}=dim (\mathfrak{g}_\lambda),$ and that $\disp \varphi_{\lambda} =
\varphi_{\mu}$ iff $\lambda = s\mu$ for some $s \in \mathfrak{w}$.  Some of
the well-known properties are $\varphi_{-\lambda}(x^{-1}) =
\varphi_{\lambda}(x)$, $\varphi_{-\lambda}(x) =
\bar{\varphi}_{\bar{\lambda}}(x),\; \lambda \in \mathfrak{a}^*_{\C},\;\;
x \in G$, and if $\Omega$ is the \textit{Casimir operator} on $G$ then
$\Omega\varphi_{\lambda} = -(\langle\lambda,\lambda\rangle +
\langle \rho, \rho\rangle)\varphi_{\lambda},$ where $\lambda \in
\mathfrak{a}^*_{\C}$ and $\langle\lambda,\mu\rangle
:=tr(adH_{\lambda} \ adH_{\mu})$ for elements $H_{\lambda}$, $H_{\mu}
\in {\mathfrak{a}}.$ The elements $H_{\lambda}$, $H_{\mu}
\in {\mathfrak{a}}$  are uniquely defined by the requirement that $\lambda
(H)=tr(adH \ adH_{\lambda})$ and $\mu
(H)=tr(adH \ adH_{\mu})$ for every $H \in {\mathfrak{a}}$ ([$6.$],
Propositions $3.1.4,\;3.2.1,\;3.2.2$ and Theorem $3.2.3$).  Clearly $\Omega\varphi_0 = 0.$

Let $$\varphi_0(x):= \int_{K}e^{(-\rho(H(xk)))}dk$$ be denoted
as $\Xi(x)$ and define $\sigma: G \rightarrow \C$ as
$\sigma(x) = \|X\|$ for every $x = k\exp X \in G,\;\; k \in K,\; X
\in \mathfrak{a},$ where $\|\cdot\|$ is a norm on the finite-dimensional
space $\mathfrak{a}.$ These two functions are zeroth elementary spherical functions on
$G$ and there exist numbers $c,d$ such that $1 \leq \Xi(a)
e^{\rho(\log a)} \leq c(1+\sigma(a))^d.$ Also there exists $r_0
> 0$ such that $c_0 =: \int_G\Xi(x)^2(1+\sigma(x))^{r_0}dx
< \infty$ ($[6.],$ p. $254$). For each
$0 \leq p \leq 2$ define ${\cal C}^p(G)$ to be the set consisting of
functions $f$ in $C^{\infty}(G)$ for which $$\mu^{p}_{g_1,
g_2;m}(f) :=\sup_G|f(g_1; x ; g_2)|\Xi (x)^{-2/p}(1+\sigma(x))^m <
\infty$$ where $g_1,g_2 \in \mathfrak{U}(\mathfrak{g}_{\C}),$ the \textit{universal
enveloping algebra} of $\mathfrak{g}_{\C},$ $m \in \Z^+, x \in G,$
$f(x;g_2) := \left.\frac{d}{dt}\right|_{t=0}f(x\cdot(\exp tg_2))$
and $f(g_1;x) :=\left.\frac{d}{dt}\right|_{t=0}f((\exp
tg_1)\cdot x).$

We call ${\cal C}^p(G)$ the Schwartz-type space on $G$
for each $0 < p \leq 2$ and note that ${\cal C}^2(G)$ is the earlier
Harish-Chandra space ${\cal C}(G)$ of rapidly decreasing functions on
$G.$ The inclusions $$C^{\infty}_{c}(G)\subset\bigcap_{0<p\leq 2}{\cal C}^p(G) \subset {\cal C}^p(G)
\subset L^p(G)$$ are continuous and with dense images. It also follows that
${\cal C}^p(G) \subseteq {\cal C}^q(G)$ whenever $0 \leq p \leq q
\leq 2.$ Each ${\cal C}^p(G)$ is closed under \textit{involution} and the
\textit{convolution}, $*.$ Indeed ${\cal C}^p(G)$ is a Fr$\acute{e}$chet algebra ($[12c.],\;p.\;357$) and the relation ${\cal C}^p(G)\ast{\cal C}^q(G)\subset{\cal C}^p(G)$ holds for all $p\geq q$ with $\frac{1}{p}+\frac{1}{q}=1;\;\;[3c.],$ Theorem $5.1.$ We endow ${\cal C}^p(G//K)$
with the relative topology as a subset of ${\cal C}^p(G).$

We shall say a function $f$ on $G$ satisfies a \textit{general strong inequality} if for any $r\geq0$ there is a constant $c=c_{r}>0$ such that
$$\mid f(y) \mid \leq c_{r} \Xi(y^{-1}x) (1+\sigma(y^{-1}x))^{-r}\;\;\;\;\;\forall\;x,y \in G.$$ We observe that if $x=e$ then, using the fact that $\Xi(y^{-1})=\Xi(y)$ and $\sigma(y^{-1})=\sigma(y),\;\forall\;y \in G,$ such a function satisfies $$\mid f(y) \mid \leq c_{r} \Xi(y^{-1}) (1+\sigma(y^{-1}))^{-r}=c_{r} \Xi(y) (1+\sigma(y))^{-r},\;\forall\;y \in G,$$ showing that a function on $G$ which satisfies a general strong inequality satisfies in particular a \textit{strong inequality} (in the classical sense of Harish-Chandra, $[12c.]$). Members of $\mathcal{C}(G)$ are those functions $f$ on $G$ for which $f(g_1; \cdot ; g_2)$ satisfies the strong inequality, for all $g_1,g_2 \in \mathfrak{U}(\mathfrak{g}_{\C}).$ We may then define $\mathcal{C}_{x}(G)$ to be those functions $f$ on $G$ for which $f(g_1; \cdot ; g_2)$ satisfies the general strong inequality, for all $g_1,g_2 \in \mathfrak{U}(\mathfrak{g}_{\C})$ and a fixed $x \in G.$ It is clear that $\mathcal{C}_{e}(G)=\mathcal{C}(G)$ and that $\bigcup_{x \in G}\mathcal{C}_{x}(G),$ which contains $\mathcal{C}(G),$ may be given an inductive limit topology.

\textbf{2.1 Proposition.} \textit{$\bigcup_{x \in G}\mathcal{C}_{x}(G)$ is a Schwartz algebra$.\Box$}

The algebra $\bigcup_{x \in G}\mathcal{C}_{x}(G)$ is worthy of an independent study. See $[9b.].$

For any measurable function $f$ on $G$ we define the \textit{Harish-Chandra Fourier
transform} $f\mapsto\mathcal{H}(f)$ as $\mathcal{H}(f)(\lambda) = \int_G f(x)
\varphi_{-\lambda}(x)dx,$ $\lambda \in \mathfrak{a}^*_{\C}=:\mathfrak{F}.$ We shall call it \textit{spherical} whenever $f\in\mathcal{C}^{p}(G//K)$ and, in this case, it may be shown that it is sufficient to define $\mathcal{H}f$ as $$(\mathcal{H}f)(\lambda):=\int_{AN}f(an)e^{(-\lambda+\rho)(\log a)}dadn,\;\lambda\in{\mathfrak{F}}_{I};$$ $[12a.],\;p.\;364.$

It is known (see $[6.]$) that for $f,g \in L^1(G)$ we have:
\begin{enumerate}
\item [$(i.)$] $\mathcal{H}(f*g) = \mathcal{H}(f)\cdot\mathcal{H}(g)$ on $ {\mathfrak{F}}^{1}$
whenever $f$ (or $g$) is right - (or left-) $K$-invariant; \item
[$(ii.)$] $\mathcal{H}(f^*)(\varphi) =
\overline{\mathcal{H}(f)(\varphi^*)},\;\varphi \in {\mathfrak{F}}^{1};$ hence
$\mathcal{H}(f^*) = \overline{\mathcal{H}(f)}$ on ${\cal P}:$ and, if we
define $f^{\#}(g) := \int_{K\times K}f(k_1xk_2)dk_1dk_2,  x\in
G,$ then \item [$(iii.)$] $\mathcal{H}(f^{\#})=\mathcal{H}(f)$ on ${\mathfrak{F}}^{1},$ where ${\mathfrak{F}}^{1}$ is the set of all bounded spherical functions and ${\cal P}$ is the subset of all positive-definite spherical functions.
\end{enumerate}

In order to know the image of the Harish-Chandra Fourier transform when
restricted to ${\cal C}^p(G//K)$ we need the following \textit{tube-spaces} that are central to the statement
of the well-known result of Trombi and Varadarajan [$11.$] (Theorem $2.2$ below).

Let $C_\rho$ be the closed convex hull of the (finite) set $\{s\rho :
s\in \mathfrak{w}\}$ in $\mathfrak{a}^*$, i.e., $$C_\rho =
\left\{\sum^n_{i=1}\lambda_i(s_i\rho) : \lambda_i \geq 0,\;\;\sum^n_{i=1}\lambda_i = 1,\;\;s_i \in \mathfrak{w}\right\}$$ where we recall that, for every
$H \in \mathfrak{a},$ $(s\rho)(H) = \frac{1}{2} \sum_{\lambda\in\triangle^+}
 m_{\lambda}\cdot\lambda (s^{-1}H).$  Now for each
$\epsilon > 0$ set ${\mathfrak{F}}^{\epsilon} = \mathfrak{a}^*+i\epsilon
C_\rho.$ Each ${\mathfrak{F}}^{\epsilon}$ is convex in $\mathfrak{a}^*_{\C}$ and
$$int({\mathfrak{F}}^{\epsilon}) =
\bigcup_{0<\epsilon'<\epsilon}{\mathfrak{F}}^{\epsilon^{'}}$$
([$11.$], Lemma $(3.2.2)$).  Let us define $\mathcal{Z}({\mathfrak{F}}^{0}) = \mathcal{S}
(\mathfrak{a}^*)$ and, for each $\epsilon>0,$ let
$\mathcal{Z}({\mathfrak{F}}^{\epsilon})$ be the space of all $\C$-valued
functions $\Phi$ such that  $(i.)$ $\Phi$ is defined and holomorphic
on $int({\mathfrak{F}}^{\epsilon}),$ and $(ii.)$ for each holomorphic
differential operator $D$ with polynomial coefficients we have $\sup_{int({\mathfrak{F}}^{\epsilon})}|D\Phi| < \infty.$ The space
$\mathcal{Z}({\mathfrak{F}}^{\epsilon})$ is converted to a Fr$\acute{e}$chet algebra by equipping it with the
topology generated by the collection, $\| \cdot \|_{\mathcal{Z}({\mathfrak{F}}^{\epsilon})},$ of seminorms given by $\|\Phi\|_{\mathcal{Z}({\mathfrak{F}}^{\epsilon})} := \sup_{int({\mathfrak{F}}^{\epsilon})}|D\Phi|.$  It is known that $D\Phi$ above extends to a continuous function on all of ${\mathfrak{F}}^{\epsilon}$
([$11.$], pp. $278-279$).  An appropriate subalgebra of
$\mathcal{Z}({\mathfrak{F}}^{\epsilon})$ for our purpose is the closed
subalgebra $\bar{\mathcal{Z}}({\mathfrak{F}}^{\epsilon})$ consisting of
$\mathfrak{w}$-invariant elements of $\mathcal{Z}({\mathfrak{F}}^{\epsilon}),$
$\epsilon \geq 0.$

{\bf 2.2 Theorem (Trombi-Varadarajan $[11.]$).}  \textit{Let $0 < p \leq 2$ and
set $\epsilon = \left(2/p\right)-1$.  Then the
Harish-Chandra Fourier transform $f \mapsto \mathcal{H}f$ is a linear
topological algebra isomorphism of ${\cal C}^p(G//K)$ onto $\bar{\mathcal{Z}}
({\mathfrak{F}}^{\epsilon}).\;\;\Box$}\\

For the Schwartz algebras ${\cal C}^p(G/K)$ a larger image than $\bar{\mathcal{Z}}({\mathfrak{F}}^{\epsilon})$ is required under the Harish-Chandra Fourier transform. Following Eguchi M. and Kowata A. $[4.]$ and Eguchi M. $[3b.]$ we define the space $\bar{\mathcal{Z}}(K/M\times{\mathfrak{F}}^{\epsilon})$ as the space of all $\mathfrak{w}$-invariant $C^{\infty}$ complex-valued functions $F$ on $K/M\times{\mathfrak{F}}_{I}$ which satisfy the following conditions:

$(i)$ for any $k\in K,$ the function $\lambda\mapsto F(kM:\lambda)$ extends holomorphically to $int({\mathfrak{F}}^{\epsilon});$

$(ii)$ for any $m\in\Z^{+},\;v \in S(\mathfrak{F}),$ $$\zeta^{\epsilon}_{v;m}(F) :=\sup_{(kM:\lambda)\in K/M\times int(\mathfrak{F})^{\epsilon}}|F(kM:\lambda ; \partial(v)|(1+ \mid \lambda\mid)^{m} <\infty.$$ The seminorms $\zeta^{\epsilon}_{v;m}$ restrict on $\bar{\mathcal{Z}}
({\mathfrak{F}}^{\epsilon})$ to the earlier Trombi-Varadarajan seminorms, $\| \cdot \|_{\mathcal{Z}({\mathfrak{F}}^{\epsilon})},$ and convert $\bar{\mathcal{Z}}(K/M\times{\mathfrak{F}}^{\epsilon})$ into a Fr$\acute{e}$chet space. Indeed, $\bar{\mathcal{Z}}({\mathfrak{F}}^{\epsilon})\subset\bar{\mathcal{Z}}(K/M\times{\mathfrak{F}}^{\epsilon}),$ as a closed subspace.

We define the map ${\cal C}^p(G/K)\rightarrow\bar{\mathcal{Z}}(K/M\times{\mathfrak{F}}^{\epsilon}):f\mapsto\mathcal{H}(f)$ now as $$\mathcal{H}(f)(kM:\lambda)=\int_{AN}f(kan)e^{(-\lambda+\rho)(\log a)}dadn,\;k\in K,\;\lambda\in\mathfrak{F}^{\epsilon},$$ referring to it as the \textit{symmetric} Harish-Chandra Fourier transform. A very important improvement on Theorem $2.2$ is the following.

{\bf 2.3 Theorem (Eguchi $[3b.]$).}  \textit{Let $0 < p \leq 2$ and
set $\epsilon = \left(2/p\right)-1$.  Then the
Harish-Chandra Fourier transform $f \mapsto \mathcal{H}f$ is a linear
topological algebra isomorphism of ${\cal C}^p(G/K)$ onto $\bar{\mathcal{Z}}(K/M\times{\mathfrak{F}}^{\epsilon}).\;\;\Box$}\\

The main result of Oyadare $[9c.]$ is a complete generalization of the above Theorems $2.2$ and $2.3$ to all $G.$ The statement is as follows. Set $$\mathcal{C}^{p}(\widehat{G}):=
\{(\mathcal{H}\xi_{1})^{-1}\cdot h\cdot(\mathcal{H}\xi_{1})^{-1}:\;h\in\bar{\mathcal{Z}}({\mathfrak{F}}^{\epsilon})\}$$ and consider $\mathcal{H}$ on all of $\mathcal{C}^{p}(G).$ We have the following.

\textbf{2.4 Theorem} (The Fundamental Theorem of Harmonic Analysis on $G$). \textit{Let $0<p\leq 2,$ then the Harish-Chandra Fourier transform, $\mathcal{H},$ sets up a linear topological algebra isomorphism $\mathcal{H}:\mathcal{C}^{p}(G)\rightarrow\mathcal{C}^{p}(\widehat{G}).$}

The above image $\mathcal{C}^{2}(\widehat{G})$ may be seen as a concrete realization of Arthur's image $\mathcal{C}^{2}(\widehat{G})_{A}$ (as computed in $[1b.,c.]$) for $\mathcal{C}^{2}(G)$ under Arthur's Harish-Chandra transform $\mathcal{H}_{A}.$ If we now set $\vartheta_{2}=\mathcal{H}_{A}\circ\mathcal{H}^{-1},$ then the properties of the topological algebra isomorphisms $\mathcal{H}_{A}$ and $\mathcal{H}$ show that $$\vartheta_{2}:\mathcal{C}^{2}(\widehat{G})\rightarrow\mathcal{C}^{2}(\widehat{G})_{A}$$ is a composition of topological algebra isomorphisms and hence that $\vartheta_{2}$ is itself a topological algebra isomorphism of our $\mathcal{C}^{2}(\widehat{G})$ with $\mathcal{C}^{2}(\widehat{G})_{A}.$ It may be of lasting importance to get explicit expressions for the map $\vartheta_{2},$ for specific examples of $G.$

The polar decomposition of $G$ implies that every $K-$biinvariant function on $G$ is completely determined by its restriction to $A^{+}.$ An example of such a function is the \textit{(zonal) spherical function,} $\varphi_{\lambda}, \lambda \in \mathfrak{a}^{\ast}_\mathbb{\C},$ on $G.$ If we denote the restriction of $\varphi_{\lambda}$ to $A^{+}$ as $\tilde \varphi_{\lambda},$ then the following system of differential equations hold: $$\tilde q \tilde \varphi_{\lambda}= \gamma (q)(\lambda) \tilde \varphi_{\lambda},$$ where $q \in \mathfrak{Q}(\mathfrak{g}_\mathbb{\C})(:= U(\mathfrak{g}_\mathbb{\C})^{K}$ = centralizer of K in $U(\mathfrak{g}_\mathbb{\C})$), $\gamma := \gamma_{\mathfrak{g}/\mathfrak{a}}$ is the \textit{Harish-Chandra homomorphism} of $\mathfrak{Q}(\mathfrak{g}_\mathbb{\C})$ onto $U(\mathfrak{g}_\mathbb{\C})^{\mathfrak{w}},$ the $\mathfrak{w}-$ invariant subspace of $U(\mathfrak{g}_\mathbb{\C}),$ with $\mathfrak{w}$ denoting the \textit{Weyl group} of the pair $(\mathfrak{g}, \mathfrak{a}),$ $\mathfrak{t} U(\mathfrak{g}_\mathbb{\C}) \bigcap \mathfrak{Q}(\mathfrak{g}_\mathbb{\C})$ is the kernel of $\gamma$ and $\tilde q$ is the restriction of $q$ to $A^{+}.$ Since $$ \widetilde{q \cdot f}
=\tilde q \cdot \tilde f,$$ for every $f \in C^{\infty}(G//K)$ we conclude that $\tilde q$ is the \textit{radial component} of $q.$ We define $q \in \mathfrak{Q}(\mathfrak{g}_\mathbb{\C})$ to be \textit{spherical} whenever $q = \tilde q.$

\ \\
\ \\
{\bf \S3. Spherical and symmetric wave-packets.}

The fundamental theorem of harmonic analysis on any locally compact group is the seeking of a linear topological algebra isomorphism of its appropriately defined Schwartz space of functions, thus characterizing its image. The implied inverse map (from the image back to the Schwartz space of the group) is called the \textit{wave-packet map,} in line with the classical case. The existence of wave-packets are always sought and are especially appropriate when the \textit{Plancherel inverse transform} is not immediate or not available. There is indeed a rich theory of wave-packets for the Schwartz-type spaces of spherical functions ($[6],$ $[10a, c]$ and $[12b.]$), of Schwartz-type functions on the symmetric spaces ($[3b.]$) and (recently) of all Schwartz-type functions on a connected semisimple Lie group ($[9c.]$). We shall review this theory in the present section.

Let $G,$ $\mathcal{C}^{p}(G//K)$ and $\bar{\mathcal{Z}}({\mathfrak{F}}^{\epsilon})$ be as in \S2. and recall the spherical (Fourier) transform of any $f\in\mathcal{C}^{p}(G//K)$ defined as $\mathcal{H}f(\lambda)=\int_{G}f(x)\varphi_{\lambda}(x)dx,$ $\lambda\in{\mathfrak{F}}^{\epsilon},$ which is a linear topological algebra isomorphism of $\mathcal{C}^{p}(G//K)$ onto $\bar{\mathcal{Z}}({\mathfrak{F}}^{\epsilon}),$ ($[6],$ p. $354$).The surjectivity part of this isomorphism (which is also called \textit{the wave-packets theorem}) requires that for any $a\in\bar{\mathcal{Z}}({\mathfrak{F}}^{\epsilon})$ the (spherical) wave-packet $\psi^{sph}_{a},$ defined as $$\psi^{sph}_{a}(x)=\frac{1}{\mid\mathfrak{w}\mid}\int_{\mathfrak{F}_{I}}
a(\lambda)\varphi_{\lambda}(x)c(\lambda)^{-1}c(-\lambda)^{-1}d\lambda$$ is in $\mathcal{C}^{p}(G//K)$ and that the map $a\mapsto\psi^{sph}_{a}$ (which is the inverse of $\mathcal{\mathcal{H}},$ up to a non-zero multiple constant) is continuous from $\bar{\mathcal{Z}}({\mathfrak{F}}^{\epsilon})$ into $\mathcal{C}^{p}(G//K).$ Eguchi $([4b.])$ has also shown that the above wave-packets theorem on $G//K$ extends to the symmetric space $G/K,$ where the symmetric wave-packets $\psi^{symm}_{a}$ are given as $$\psi^{symm}_{a}(kM;x)=\frac{1}{\mid\mathfrak{w}\mid}\int_{\mathfrak{F}_{I}}
a(kM;\lambda)\varphi_{\lambda}(x)c(\lambda)^{-1}c(-\lambda)^{-1}d\lambda,$$ with $a\in\bar{\mathcal{Z}}(K/M\times{\mathfrak{F}}^{\epsilon}).$ Oyadare $([9c.])$ recently gives the proof of the (full) fundamental theorem on $G,$  where the full wave-packets theorem was proved on all of $G$ for the general wave-packets, $\psi_{a},$ that are given as $$\psi_{a}=\mathcal{H}^{-1}((\mathcal{H}\xi_{1})^{-1})\ast\mathcal{H}^{-1}((\mathcal{H}\xi_{1})\cdot a\cdot(\mathcal{H}\xi_{1}))\ast\mathcal{H}^{-1}((\mathcal{H}\xi_{1})^{-1})$$ with $a\in\mathcal{C}^{p}(\widehat{G}).$
\ \\
{\bf \S4. Canonical wave-packets on $G.$}

We shall refer to the above general wave-packets, $\psi_{a},$ as \textit{the canonical wave-packets on $G.$}

{\bf 4.1 Lemma.} \textit{Let $a\in\mathcal{C}^{p}(\widehat{G}),$ $0\leq p\leq2.$ The canonical wave-packets $\psi_{a}$ on $G$ is given by the convolutions of a symmetric, a spherical and a symmetric wave-packets.}

{\bf Proof.} It is clear from the above that $\mathcal{H}^{-1}((\mathcal{H}\xi_{1})^{-1})\in\mathcal{C}^{p}(G/K)$ and that $\mathcal{H}^{-1}((\mathcal{H}\xi_{1})\cdot a\cdot(\mathcal{H}\xi_{1}))\in\mathcal{C}^{p}(G//K).\;\Box$

We note here that the symmetric part, $\mathcal{H}^{-1}((\mathcal{H}\xi_{1})^{-1}),$ of the canonical wave-packets is a fixed function in $\mathcal{C}^{p}(G/K).$ We may then have the following.

{\bf 4.2 Theorem.} \textit{Let $0\leq p\leq2,$ then $$\mathcal{C}^{p}(G)\cong\mathcal{H}^{-1}((\mathcal{H}\xi_{1})^{-1})
\ast\mathcal{C}^{p}(G//K)\ast\mathcal{H}^{-1}((\mathcal{H}\xi_{1})^{-1}).$$}

{\bf Proof.} We already know from Trombi-Varadarajan theorem that $$\mathcal{H}(\mathcal{C}^{p}(G//K))\cong\bar{\mathcal{Z}}({\mathfrak{F}}^{\epsilon}).\;\;\mbox{That is,}\;\; \mathcal{C}^{p}(G//K)\cong\mathcal{H}^{-1}(\bar{\mathcal{Z}}({\mathfrak{F}}^{\epsilon})).$$ Hence, by the continuity of the convolution, we have $$\mathcal{H}^{-1}((\mathcal{H}\xi_{1})^{-1})
\ast\mathcal{C}^{p}(G//K)\ast\mathcal{H}^{-1}((\mathcal{H}\xi_{1})^{-1})$$ $$\cong\mathcal{H}^{-1}((\mathcal{H}\xi_{1})^{-1})
\ast\mathcal{H}^{-1}(\bar{\mathcal{Z}}({\mathfrak{F}}^{\epsilon}))\ast\mathcal{H}^{-1}((\mathcal{H}\xi_{1})^{-1})$$
$$\cong\mathcal{H}^{-1}(((\mathcal{H}\xi_{1})^{-1})
\cdot(\bar{\mathcal{Z}}({\mathfrak{F}}^{\epsilon}))\cdot((\mathcal{H}\xi_{1})^{-1}))
\cong\mathcal{H}^{-1}(\mathcal{C}^{p}(\widehat{G}))\cong\mathcal{C}^{p}(G).\;\Box$$

The last theorem gives a decomposition of the Harish-Chandra Schwartz-type algebras $\mathcal{C}^{p}(G)$ showing, once again, the centrality of the spherical part $\mathcal{C}^{p}(G//K).$

It now follows that the canonical wave-packets $\psi_{a}$ on $G$ may be seen as either $$\psi_{a}=\mathcal{H}^{-1}((\mathcal{H}\xi_{1})^{-1})\ast\psi^{sph}_{a}\ast\mathcal{H}^{-1}((\mathcal{H}\xi_{1})^{-1}),$$ for any $a\in\bar{\mathcal{Z}}({\mathfrak{F}}^{\epsilon})$ or as $$\psi_{a}=\psi^{symm}_{a}\ast\mathcal{H}^{-1}((\mathcal{H}\xi_{1})^{-1}),$$ for any $a\in\bar{\mathcal{Z}}(K/M\times{\mathfrak{F}}^{\epsilon}).$ Clearly $\psi^{symm}_{a}=\mathcal{H}^{-1}((\mathcal{H}\xi_{1})^{-1})\ast\psi^{sph}_{a}$ with $\psi^{sph}_{a}$ given as $\psi^{sph}_{a}=\mathcal{H}^{-1}((\mathcal{H}\xi_{1})\cdot a\cdot(\mathcal{H}\xi_{1}))=\xi_{1}\ast \mathcal{H}^{-1}a\ast\xi_{1}.$

The above formula $\psi_{a}=\psi^{symm}_{a}\ast\mathcal{H}^{-1}((\mathcal{H}\xi_{1})^{-1})$ therefore gives the canonical wave-packets $\psi_{a}$ on $G$ as the convolution of two integrals, $\psi^{symm}_{a}$ and $\mathcal{H}^{-1}((\mathcal{H}\xi_{1})^{-1}).$ A way forward is to have a complete explicit integral expression for $\mathcal{H}^{-1}((\mathcal{H}\xi_{1})^{-1})$ and be able to compute the said convolution of the two integrals. Another approach is to study $\psi_{a}$ directly from the properties of $\psi^{symm}_{a}$ and $\mathcal{H}^{-1}((\mathcal{H}\xi_{1})^{-1})$ in the defining convolution. Indeed, for any $x\in G,$
$\psi_{a}(x)=\int\int_{G}\psi^{symm}_{a}(y)\mathcal{H}^{-1}((\mathcal{H}\xi_{1})^{-1})(y^{-1}x)dy,$ where $dy$ is a normalized Haar measure on $G.$ An explicit integral realization of the defining convolutions of the canonical wave-packets is however more desirous in order to calculate its measure, which will be the generalization of the spherical Plancherel measure $\mid c(\lambda)\mid^{-2}d\lambda$ to all of $\mathcal{C}^{2}(G)$ (or even $\mathcal{C}^{p}(G)$), leading directly to the complete form of the $c-$function.

We now have the following result.

{\bf 4.3 Proposition.} \textit{The canonical wave-packets $\psi_{a}$ transform under the operators $q_{1},q_{2}\in\mathcal{U}(\mathfrak{g}_{\C})$ as $$q_{1}\psi_{a}q_{2}=\psi^{symm}_{a}q_{2}\ast q_{1}\mathcal{H}^{-1}((\mathcal{H}\xi_{1})^{-1}).$$ In particular, $\psi_{a}q=\gamma(q)(\lambda)\psi_{a},$ for every $q\in\mathcal{U}(\mathfrak{g}_{\C})^{K}.$}

{\bf Proof.} Clear from $(6.1.9)$ of $[6].\;\Box$

We may now consider the wave-packet map $W$ on $\mathcal{C}^{2}(\widehat{G})$ given as $$W:a\mapsto\psi^{symm}_{a}\ast \mathcal{H}^{-1}((\mathcal{H}\xi_{1})^{-1}).$$ Set $\mathcal{C}^{\#}(G):=\{Wa:a\in\mathcal{C}^{2}(\widehat{G})\}$ which is the image of $W$ in $\mathcal{C}^{2}(G)$ and $^{0}\mathcal{C}(G):=kernel(W\mathcal{H})$ in $\mathcal{C}^{2}(G).$

{\bf 4.4 Theorem.} \textit{$^{0}\mathcal{C}(G)$ is the orthogonal complement of $\mathcal{C}^{\#}(G)$ in $\mathcal{C}^{2}(G).$ In particular, $\mathcal{C}^{2}(G)=^{0}\mathcal{C}(G)\oplus\mathcal{C}^{\#}(G),$ with continuous projections.}

{\bf Proof.} As in Theorem $60$ of $[12c.].\;\Box$

The topological algebra isomorphism $\mathcal{H}:\mathcal{C}^{2}(G)\cong\mathcal{C}^{2}(\widehat{G})$ assures us that the above decomposition in Theorem $4.4$ also transfers to $\mathcal{C}^{2}(\widehat{G})$ to give $\mathcal{C}^{2}(\widehat{G})=^{0}\mathcal{C}(\widehat{G})\oplus\mathcal{C}^{\#}(\widehat{G}),$ with continuous projections.
\ \\
{\bf   References.}
\begin{description}
\item [{[1.]}] Arthur, J. G., $(a.)$ \textit{Harmonic analysis of tempered distributions on semi-simple Lie groups of real rank one}, Ph.D. Dissertation, Yale University, $1970;$ $(b.)$ \textit{Harmonic analysis of the Schwartz space of a reductive Lie group I,} mimeographed note, Yale University, Mathematics Department, New Haven, Conn; $(c.)$ \textit{Harmonic analysis of the Schwartz space of a reductive Lie group II,} mimeographed note, Yale University, Mathematics Department, New Haven, Conn.
        \item [{[2.]}] Barker, W. H.,  $L^p$ harmonic analysis on $SL(2, \R),$ \textit{Memoirs of American Mathematical Society,} $vol.$ {\bf 76} , no.: {\bf 393}. $1988$
            \item [{[3.]}] Eguchi, M., $(a.)$ The Fourier Transform of the Schwartz space on a semisimple Lie group, \textit{Hiroshima Math. J.}, $vol.$ \textbf{4}, ($1974$), pp. $133-209.$ $(b.)$ Asymptotic expansions of Eisenstein integrals and Fourier transforms on symmetric spaces, \textit{J. Funct. Anal.} \textbf{34}, ($1979$), pp. $167 - 216.$ $(c.)$ Some properties of Fourier transform on Riemannian symmetric spaces, \textit{Lecture on Harmonic Analysis on Lie Groups and related Topics,} (T. Hirai and G. Schiffmann (eds.)) Lectures in Mathematics, Kyoto University, No. \textbf{4}) pp. $9 - 43.$
                \item [{[4.]}] Eguchi, M. and Kowata, A., On the Fourier transform of rapidly decreasing function of $L^{p}$ type on a symmetric space, \textit{Hiroshima Math. J.} $vol.$ \textbf{6}, ($1976$), pp. $143 - 158.$
                \item [{[5.]}] Ehrenpreis, L. and Mautner, F. I., Some properties of the Fourier transform on semisimple Lie groups, I,
                \textit{Ann. Math.}, $vol.$ $\textbf{61}$ ($1955$), pp. 406-439; II, \textit{Trans. Amer. Math. Soc.}, $vol.$ $\textbf{84}$ ($1957$), pp. $1-55;$ III, \textit{Trans. Amer. Math. Soc.}, $vol.$ $\textbf{90}$ ($1959$), pp. $431-484.$
                    \item [{[6.]}] Gangolli, R. and Varadarajan, V. S., \textit{Harmonic analysis of spherical functions on real reductive groups,} Ergebnisse der Mathematik und iher Genzgebiete, $vol.$ {\bf 101}, Springer-Verlag, Berlin-Heidelberg. $1988.$
                    \item [{[7.]}] Helgason, S., $(a)$  A duality for symmetric spaces with applications to group representations, \textit{Advances in Mathematics,} $vol.$ $\textbf{5}$ ($1970$), pp. $1-154.$ $(b)$ \textit{Differential geometry and symmetric spaces,} Academic Press, New York, $1962.$
                        \item [{[8.]}] Knapp, A.W., \textit{Representation theory of semisimple groups; An overview based on examples,} Princeton University Press, Princeton, New Jersey. $1986.$
                \item [{[9.]}] Oyadare, O. O., $(a.)$ On harmonic analysis of spherical convolutions on semisimple Lie groups, \textit{Theoretical Mathematics and Applications}, $vol.$ $\textbf{5},$ no.: {\bf 3}. ($2015$), pp. $19-36.$ $(b.)$ Series analysis and Schwartz algebras of spherical convolutions on semisimple Lie groups (\textit{under review}), See arXiv.$1706.09045$ [math.RT]. $(c.)$ Non-spherical Harish-Chandra Fourier transforms on real reductive groups (\textit{under review}), See arXiv.$1907.00717$ [math.FA]. $(d.)$ The full Bochner theorem on real reductive groups (\textit{under review}), See arXiv.$1907.10819$ [math.FA].
                    \item [{[10.]}] Trombi, P. C., $(a.)$ Spherical transforms on symmetric spaces of rank one (or Fourier analysis on semisimple Lie groups of split rank one), \textit{Thesis, University of Illinios} ($1970$). $(b.)$ On Harish-Chandra's theory of the Eisenstein integral for real semisimple Lie groups. \textit{University of Chicago Lecture Notes in Representation Theory,} ($1978$), pp. $287$-$350.$ $(c.)$ Harmonic analysis of $\mathcal{C}^{p}(G:F)\;(1\leq p<2)$ \textit{J. Funct. Anal.,} $vol.$ {\bf 40}. ($1981$), pp. $84$-$125.$ $(d.)$ Invariant harmonic analysis on split rank one groups with applications. \textit{Pacific J. Math.} $vol.$ \textbf{101}. no.: \textbf{1}.($1982$), pp. $223 - 245.$
                \item [{[11.]}] Trombi, P. C. and Varadarajan, V. S., Spherical transforms on semisimple Lie groups, \textit{Ann. Math.,} $vol.$ {\bf 94}. ($1971$), pp. $246$-$303.$
                        \item [{[12.]}] Varadarajan, V. S., $(a.)$ Eigenfunction expansions on semisimple Lie groups, in \textit{Harmonic Analysis and Group Representation}, (A. Fig$\grave{a}$  Talamanca (ed.)) (Lectures given at the $1980$ Summer School of the \textit{Centro Internazionale Matematico Estivo (CIME)} Cortona (Arezzo), Italy, June $24$ - July $9.$ vol. \textbf{82}) Springer-Verlag, Berlin-Heidelberg. $2010,$ pp. $351-422.$ $(b.)$ \textit{An introduction to harmonic analysis on semisimple Lie groups,} Cambridge Studies in Advanced Mathematics, \textbf{161},  Cambridge University Press, $1989.$ $(c.)$ Harmonic analysis on real reductive reductive groups, \textit{Lecture Notes in Mathematics,} \textbf{576}, Springer Verlag, $1977.$
                        \end{description}
\end{document}